\documentclass[reqno]{gokova}
\usepackage{epsfig,graphicx}

\newtheorem{thm}{Theorem}[section]

\newtheorem{lem}[thm]{Lemma}
\newtheorem{prop}[thm]{Proposition}

\theoremstyle{remark}
\newtheorem{rem}{\rm\bfseries{Remark}}[section]

\DeclareMathSymbol{\pitchfork}    {\mathrel}{AMSa}{"74}
\DeclareMathSymbol{\varsubsetneq}   {\mathrel}{AMSb}{"20}
\DeclareMathSymbol{\smallsetminus}  {\mathbin}{AMSb}{"72}
\def\N{\ifmmode{\mathbb N}\else{$\mathbb N$}\fi} 
\def\R{\ifmmode{\mathbb R}\else{$\mathbb R$}\fi} 
\def\Q{\ifmmode{\mathbb Q}\else{$\mathbb Q$}\fi} 
\def\C{\ifmmode{\mathbb C}\else{$\mathbb C$}\fi} 
\def\Z{\ifmmode{\mathbb Z}\else{$\mathbb Z$}\fi} 
\def\la{\langle}
\def\ra{\rangle}
\def\neqqq{\varsubsetneq}
\def\nologo{\let\logo@\relax}
\def\stmin{\smallsetminus}

\def\d{\partial}
\def\:{:}
\def\conj{\text{\rm{conj}}}
\def\ge{\geq}
\def\le{\leq}
\def\Cl{\rm{Cl}}
\def\Int{\rm{Int}}
\def\trans{\pitchfork}
\def\vbot{\bot\!\!\!\bot}
\def\Cp#1{\mathbb{C}\rm{P}^{#1}}
\def\Rp#1{\mathbb{R}\rm{P}^{#1}}
\def\barCP#1{\overline{\C\rm P}^{#1}}
\def\b{\beta}
\def\D{\Delta}
\newcommand{\e}{\varepsilon}
\def\G{\Gamma}

\begin{document}

\title{Knotting of algebraic curves in complex surfaces}
\author[S. Finashin \ \  Knotting of Algebraic Curves]{Sergey Finashin}
\begin{abstract}
For any $d\ge 5$, I constructed infinitely many pairwise smoothly
non-equivalent surfaces $F\subset\Cp{2}$ homeomorphic to a non-singular
algebraic curve of degree $d$, realizing the same homology class as such
a curve and having  abelian fundamental group $\pi_1(\Cp2\stmin F)$.

It is a special case of a more general theorem, which concerns for instance
those algebraic curves, $A$,
in a simply connected algebraic surface, $X$,
which admit irreducible degenerations to a curve
$A_0$, with a unique singularity of the type
$X_9$, and such that  $A\circ A>16$.
\end{abstract}
\address{Middle East Technical University, Ankara 06531 Turkey}
\email{serge@metu.edu.tr}
%
%
\volume{7}
\maketitle
\section{Introduction}

\begin{thm}\label{thm-for-plane}
For any $d\ge5$ there exist infinitely many
smooth oriented closed surfaces
$F\subset\Cp2$ representing class $d\in H_2(\Cp2)=\Z$,
having $\rm{genus}(F)=\frac12(d-1)(d-2)$
and $\pi_1(\Cp2\stmin F)\cong\Z/d$, such that
the pairs $(\Cp2,F)$ are pairwise smoothly non-equivalent.
Moreover, $d$-fold cyclic coverings over $\Cp2$ branched along
$F$ differ by their Seiberg-Witten invariants and thus are non-diffeomorphic.
\end{thm}


This theorem, which answers an old question (cf. \cite{K}, Problem 4.110),
is proved in \cite{F2} for even $d\ge6$.
In this paper I added the proof for odd $d$ and generalized Theorem
\ref{thm-for-plane} (see below Theorem \ref{general}).
Sections 2-3 and the Appendix reproduce the content of \cite{F2}
whereas Section 5 extends the results from there.

\begin{rem}
Note that the surfaces that I construct are not symplectic.
Some speculation referring to Gromov's theorem suggests that
any symplectic surface in $\Cp2$ may be isotopic to an algebraic curve.
As far as I know, at the moment it is proved only for degrees $d\le4$.
\end{rem}

The knotting construction used to obtain surfaces $F$ is a
relative of the rim-surgery defined in \cite{FS2}. An alternative
way to achieve Theorem \ref{thm-for-plane} is to use the
tangle-surgery of Viro introduced in \cite{FKV}. For technical
reasons I prefer to use the rim-surgery in this paper, and give below
an idea about the other approach just because it inspired this paper.

\subsection{The idea that inspired my construction}
Any kind of a surgery on a codimension two submanifold, $F$, in
some fixed $n$-manifold $X$ gives rise to some $n$-dimensional
surgery on the double covering $Y\to X$ branched along $F$. Vice
versa, considering a surgery on $Y$, one can try to perform it
equivariantly with respect to the covering transformation, which
results in some surgery on a pair $(X,F)$. Sometimes $X$ is
preserved, and only $F$ as an embedded submanifold is modified by
this surgery. I call such an ambient surgery on $F$ in $X$ {\it
the folding} of the corresponding surgery on $Y$.

For example, if $Y$ is a complex surface defined over $\R$, and
$X=Y/\conj$ is the quotient by the complex conjugation $\conj\:
Y\to Y$, then the projection $p\: Y\to X$ is a double covering
branched along $F=\rm{Fix}(\conj)$  (the real locus of $Y$).
Algebraic transformations (say, a blow-up, or a logarithmic
transform) can be applied to $Y$ in the real category. It turns
out (at least in the examples known to the author) that the
quotient $X=Y/\conj$ is not changed if a transformation is
irreducible over $\C$, .i.e., if it does not contain a pair of
$\conj$-symmetric transformations localized outside the real part
$F$.

Say, the folding of a blow-up at a real point of $Y$ is a real
blow-up of $F$, that is an ambient connected sum $(X,F)\#(S^4,\Rp2)$,
because $\Cp2/\conj\cong S^4$.
Viro observed \cite{FKV} that the folding of a logarithmic transform
is a certain tangle-surgery on $F$. This yields
``exotic knottings'' of $F=\#_{10}\Rp2$ in $S^4=Y/\conj$, where
$Y=E(1)=\Cp2\#_9\barCP2$ is a rational elliptic surface,
being modified by logarithmic transforms
(which produce Dolgachev surfaces defined over $\R$).

The same construction applied to a K3 surface, $Y=E(2)$, instead
of $E(1)$, gives ``exotic knottings'' of $F=\rm{Fix}(\conj)$ in
$X=Y/\conj$. For a suitable choice of the real structure in $Y$,
the quotient $X$ is diffeomorphic to $\Cp2$ and $F$ becomes a
sextic in $X$, so the surgery gives examples for $d=6$ in Theorem
\ref{thm-for-plane}.
Viro's tangle surgery can be applied, in general, along any
null-framed annulus membrane on a surface in a four-manifold,
which gives in the covering space a logarithmic transform.
Suitable membranes on algebraic curves in $\Cp2$ are described in
what follows.

It turned out that the Fintushel-Stern's surgery on $Y$ admits
also a folding, i.e., can be made equivariantly, with the
quotient $X$ being preserved, provided the knot that we use is a
double knot, i.e., $K\#K$. This folding is just what I call below
``an annulus rim surgery''.

\subsection{An annulus rim-surgery}
\begin{figure}[htb]
\includegraphics{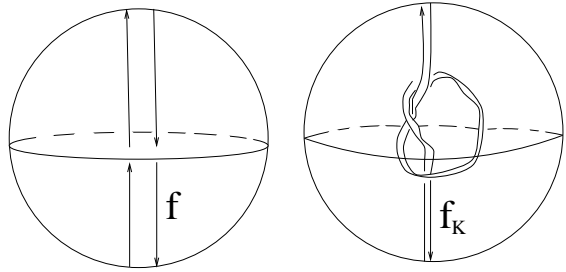}
\caption{Knotting of a band $b_K$}
\label{knotting}
\end{figure}

Our surgery, like the Viro tangle surgery, requires a suitable
annulus membrane and produces a new surface via knotting an old
one along such a membrane. By an annulus membrane for a smooth
surface $F$ in a $4$-manifold $X$ I mean a smoothly embedded
surface $M\subset X$, $M\cong S^1\times I$, with $M\cap F=\d M$
and such that $M$ comes to $F$ normally along $\d M$.
 Assume that such a membrane has
framing $0$, or equivalently, admits a diffeomorphism
of its regular neighborhood $\phi\:U\to S^1\times D^3$
 mapping $U\cap F$ onto $S^1\times f$,
where $f=I\vbot I\subset D^3$ is a disjoint union of two segments,
which are unknotted and unlinked in $D^3$,
that is to say that
 a union of $f$ with a pair of arcs on a sphere $\d D^3$
bounds a trivially
embedded band, $b\subset D^3$, $b\cong I\times I$,  so that
$f= I\times(\d I)\subset b$ (see Figure \ref{knotting}).
The annulus $M$ can be viewed
as $S^1\times\{\frac12\}\times I$ in
$S^1\times b\subset S^1\times D^3\cong U$.

If $X$ and $F$ are oriented, then $f$ inherits an orientation
as a transverse intersection, $f=F\trans D^3$, and
we may choose a band $b$ so that the orientation of $f$
is induced from some orientation of $b$.
It is convenient to view $f=I\vbot I$ as is shown on Figure \ref{knotting},
so that the segments of $f$ are parallel and oppositely oriented,
with $b$ being a thin band between them.
Such a presentation is always possible if we allow a modification
of $\phi$, since one of the segments of $f$ may be turned around
by a diffeomorphism of $D^3\to D^3$ leaving the other segment fixed.

 Given a knot $K\subset S^3$, we construct a new smooth surface,
 $F_{K,\phi}$, obtained from $F$ by tying a pair of segments
$I\vbot I$ along $K$ inside $D^3$, as is shown on Figure \ref{knotting}.
More precisely, we consider a band $b_K\subset D^3$ obtained from
$b$ by knotting along $K$
and let $f_K$ denote the pair of arcs bounding $b_K$ inside $D^3$.
We assume that the framing of
$b_K$ is chosen the same as the framing of $b$, or equivalently,
that the inclusion homomorphisms from
$H_1(\d D^3\stmin (\d f))=H_1(\d D^3\stmin (\d f_K))$
to $H_1(B^3\stmin f)$ and to $H_1(B^3\stmin f_K)$ have the same kernel.
Then $F_{K,\phi}$ is obtained from $F$ by replacing
$S^1\times f\subset S^1\times D^3\cong U$ with
$S^1\times f_K$.
It is obvious that
$F_{K,\phi}$ is homeomorphic to $F$ and realizes
the same homology class in $H_2(X)$.

The above construction is called in what follows
{\it an annulus rim-surgery}, since it looks like the
 rim-surgery of Fintushel and Stern \cite{FS2}, except that
we tie two strands simultaneously, rather then one.
Recall that the usual rim-surgery is applied in \cite{FS2} to
surfaces $F\subset X$ which are {\it primitively embedded},
 that is $\pi_1(X\stmin F)=0$, which is not the case for
the algebraic curves in $\Cp2$ of degree $>1$.
The primitivity condition  is required to
preserve the fundamental group of $X\stmin F$ throughout the knotting.
 An annulus rim-surgery may preserve {\it a non-trivial}
group $\pi_1(X\stmin F)$, if we require
commutativity of $\pi_1(X\stmin(F\cup M))$, instead of primitivity
of the embedding.

\begin{prop}\label{1.2}
Assume that $X$ is a simply connected closed $4$-manifold,
$F\subset X$ is an oriented closed surface with an annulus-membrane
$M$ of index $0$, $\phi\: U\to S^1\times D^3$ is a trivialization like
described above and $K\subset S^3$ is any knot.
Assume furthermore that
 $F\stmin\d M$ is connected and the group $\pi_1(X\stmin(F\cup M))$
is abelian.
Then the group $\pi_1(X\stmin F_{K,\phi})$ is cyclic and isomorphic to
$\pi_1(X\stmin F)$.
\end{prop}

\subsection{Maximal nest curves}
To prove Theorem \ref{thm-for-plane}, I apply an annulus rim-surgery
inside $X=\Cp2$ letting $F=\C A$ be the complex point set of
a suitable non-singular real algebraic curve,
containing an annulus, $M$, among
the connected components of $\Rp2\stmin\R A$,
where $\R A=\C A\cap\Rp2$ is the real locus of the curve.

One may take, for instance, a real algebraic curve
$\C A$ of degree $d$, with {\it a maximal nest real scheme}.
 Such a curve for $d=2k$ is constructed by a small real perturbation of
a union of $k$ real conics, whose real parts (ellipses) are
ordered by inclusion in $\Rp2$.
 For $d=2k+1$, we add to such conics a real line not intersecting
the conics in $\Rp2$ and then perturb the unions. The real part,
$\R A$, of our non-singular curve contains $k$ components,
$O_1,\dots,O_k$, called ovals (just deformed ellipses). We order
the ovals so that $O_{i}$ lies inside $O_{i+1}$ and denote by
$R_i$ the annulus-component of $\Rp2\stmin\R A$ between $O_{i}$
and $O_{i+1}$ for $i=1,\dots,k-1$. $R_0$ is a topological disk
bounded from outside by $O_1$, and $R_{k}$ is the component
bounded from inside by $O_k$.

The closures, $\Cl(R_i)$, for $i=1,\dots,k-1$
are obviously $0$-framed annulus-membranes
on $\C A$. For simplicity, let us choose $M=\Cl(R_1)$.

\begin{prop}\label{1.3}
The assumptions of Proposition \ref{1.2}
 hold if we put $X=\Cp2$, let
$F=\C A$ be a maximal nest real algebraic curve of degree
 $d\ge5$ and choose $M=\Cl(R_1)$.
\end{prop}

\subsection
{ Proof of Theorem \ref{thm-for-plane} for even $d$}
Assuming that the class $[F]\in H_2(X;\Z/2)$ vanishes,
one can consider a double covering $p\:Y\to X$ branched along $F$;
such a covering is unique if we require in addition that $H_1(X;\Z/2)=0$.
 Similarly, we consider the double coverings
$Y(K,\phi)\to X$ branched along $F_{K,\phi}$.
To prove non-equivalence of pairs $(\Cp2,F_{K,\phi})$
for some family of knots $K$,
it is enough to show that $Y(K,\phi)$ are not pairwise diffeomorphic.
To show it,
I use that $Y(K,\phi)$ is diffeomorphic to the $4$-manifolds
$Y_{K\#K}$ obtained from $Y$
by a surgery introduced in \cite{FS1}
(I call it {\it FS-surgery}).

\begin{prop}\label{1.4}
The above $Y(K,\phi)$ is diffeomorphic to
a  $4$-manifold obtained from $Y$ by the FS-surgery along
the torus $T=p^{-1}(M)$ via the knot $K\# K\subset S^3$.
\end{prop}

To distinguish the diffeomorphism types of $Y_{K\#K}$
one can use the
 formula of Fintushel and Stern \cite{FS1} for SW-invariants of
a $4$-manifold $Y$ after FS-surgery along a torus $T\subset Y$.
Recall that this formula can be applied if the
SW-invariants of $Y$ are well-defined and a torus $T$,
realizing a non-trivial class $[T]\in H_2(Y)$, is {\it c-embedded}
(the latter means that $T$ lies as a non-singular fiber
in a cusp-neighborhood in $Y$, cf. \cite{FS1}).
Being an algebraic surface of genus $\ge1$, the double plane
$Y$ has well-defined
SW-invariants. The conditions on $T$ are also satisfied.

\begin{prop}\label{c-embedded}
Assume that $X$, $F$ and $M$ are like in Proposition \ref{1.2},
$[F]\in H_2(X;\Z/2)$ vanishes and $p\:Y\to X$ is like above.
Then the torus $T=p^{-1}(M)$ is primitively embedded in $Y$
and therefore $[T]\in H_2(Y)$ is an infinite order class.
If, moreover, $X$, $F$ and $M$ are chosen like in Proposition \ref{1.3}, then
$T\subset Y$ is c-embedded.
\end{prop}

Recall that the product formula \cite{FS1}
$$SW_{Y_K}=SW_Y{}^{.}\D_K(t),\ \text{ where } t=\exp(2[T])
$$
expresses the Seiberg-Witten invariants
(combined in a single polynomial) of the manifold $Y_K$, obtained
by an FS-surgery, in terms of the Seiberg-Witten invariants of $Y$
and the Alexander polynomial, $\D_K(t)$, of $K$.

 This formula implies that the basic classes of $Y_K$ can be expressed
as $\pm \b+ 2n[T]$, where
$\pm\b\in H_2(Y)$ are the basic classes of $Y$ and
 $|n|\le\deg(\D_K(t))$,
are the degrees of the non-vanishing monomials in $\D_K(t)$.
So, if $[T]$ has infinite order,
then the manifolds $Y(K,\phi)\cong Y_{K\#K}$ differ from each other by
their SW-invariants,
and moreover, by the numbers of their basic classes,
for an infinite family of knots $K$,
since the number of the basic classes is determined by
the number of the terms in $\D_{K\# K}=(\D_K)^2$
(one can take any family of knots with Alexander polynomials of
distinct degrees).
\qed

\subsection
{ A generalization}
More generally, one can produce ``fake algebraic curves'' under
the following conditions.

\begin{thm}\label{general}
Assume that $F$ is a non-singular connected curve in a
simply connected complex surface $X$, which admits a deformation
degenerating $F$ into an irreducible curve $F_0\subset X$,
with a singularity of the type $X_9$, such that
the fundamental group $\pi_1(X\stmin F_0)$ is abelian.
Then there exists an infinite family of surfaces
$F_{K,\phi}\subset X$ homeomorphic
to $F$ and realizing the same homology class as $F$, having the same
fundamental group of the complement, but with the smoothly non-equivalent
pairs $(X,F_{K,\phi})$.
\end{thm}

I remind that $X_9$-singularity is a point where
$4$ non-singular branches meet pairwise transversally.
Nori's theorem \cite{N} gives conditions under which
$\pi_1(X\stmin F_0)$ must be abelian. For instance, it is so if
$A_0$ has no other singularities except $X_9$ and $A\circ A>16$.

\begin{rem}
The claim of Theorem \ref{general}
holds also if $F_0$ has a more complicated then $X_9$ singularity,
provided the group $\pi_1(X\stmin F_0)$ is abelian.
\end{rem}

\section{Commutativity of the fundamental group throughout the knotting}

\begin{lem}\label{fg-cyclic}
The assumptions of Proposition \ref{1.2} imply that
$\pi_1(X\stmin(F\cup M))=\pi_1(X\stmin F)$
is cyclic with a generator presented
by a loop around $F$.
\end{lem}

\begin{proof}
The Alexander duality in $X$ combined with the exact cohomology
sequence of a pair $(X,F\cup M)$ gives
$$
H_1(X\stmin(F\cup M))\cong H^3(X,F\cup M)=H^2(F\cup M)/i^*
H^2(X)
$$
where $i\: F\cup M\to X$ is the inclusion map.
 If $F$ is oriented and $F\stmin\d M$ is connected, then
the Mayer-Vietoris Theorem yields
$H^2(F\cup M)\cong H^2(F)\cong\Z$, and thus
$H_1(X\stmin(F\cup M))\cong H_1(X\stmin F)$ is cyclic
with a generator presented by a loop around $F$.
The same property holds for  the fundamental groups of
$X\stmin(F\cup M)$ and $X\stmin F$,
since they are abelian by the assumption of Proposition \ref{1.2}.
\end{proof}

\noindent{\bf Proof of Proposition \ref{1.2}.}
Put $X_0=\Cl(X\stmin U)$. Then $\d X_0=\d U\cong S^1\times S^2$ and
$\d U\stmin F$ is a deformational retract of
$U\stmin(F\cup M)$, so
$$
\pi_1(X_0\stmin F)=\pi_1(X\stmin(F\cup M))
$$
Since this group is cyclic and is generated by a loop around $F$,
the inclusion homomorphism
$h\:\pi_1(\d U\stmin F)\to\pi_1(X_0\stmin F)$
is epimorphic and thus
$\pi_1(X_0\stmin F)=\pi_1(\d U\stmin F)/k$, where $k$
is the kernel of $h$.

\begin{figure}[htb]
\includegraphics{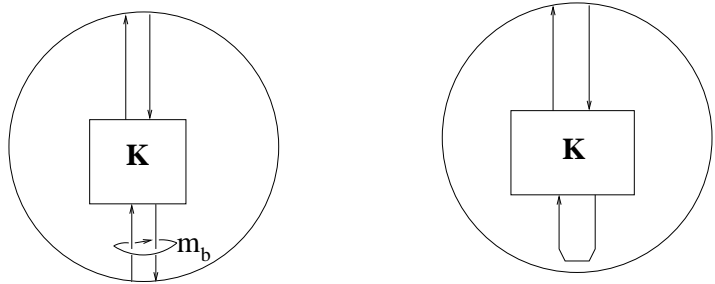}
\caption{Gluing a $2$-cell along $m_b$ effects as transforming
 $f_K$ into an unknotted arc}
\label{unknotting}
\end{figure}

Applying the Van Kampen theorem to the triad
$(X_0\stmin F,U\stmin F_{K,\phi},\d U\stmin F)$,
we conclude that
$
\pi_1(X\stmin F_{K,\phi})\cong\pi_1(U\stmin F_{K,\phi})/j(k)
$,
where $j\:\pi_1(\d U\stmin F)\to\pi_1(U\stmin  F_{K,\phi})$
is the inclusion homomorphism.
Furthermore, in the splitting
$$
\pi_1(U\stmin  F_{K,\phi})\cong\pi_1(S^1\times(D^3\stmin
f_K))\cong\Z\times\pi_1(D^3\stmin f_K)
$$
factorization by $j(k)$ kills the first factor $\Z$ and adds some
relations to $\pi_1(D^3\stmin f_K)$, one of which
effects to $\pi_1(D^3\stmin f_K)$ as if we attach
a $2$-cell  along a loop, $m_b$, going once
around the band $b_K$
(to see it, note that factorization by $k$ leaves only one generator of
$\pi_1(\d D^3\stmin f_K)=\pi_1(S^2\stmin\{4 \text{\rm{pts}}\})$).
Attaching such a $2$-cell effects to $\pi_1$
as connecting together a pair of the endpoints
of $f_K$, which transforms $f_K$ into an arc (see Figure \ref{unknotting}).
This arc is unknotted and thus factorization by $j(k)$ makes
$\pi_1(D^3\stmin f_K)$ cyclic and leaves
$\pi_1(X\stmin F_{K,\phi})$ isomorphic to
$
\pi_1(X_0\stmin F)\cong\pi_1(X\stmin(F\cup
M))\cong\pi_1(X\stmin F)
$.
\qed

\noindent{\bf Proof of Proposition \ref{1.3}.}
All the assumptions of Proposition \ref{1.2} except the last two
are obviously satisfied. It is well known that $\C A\stmin\R A$
splits for a maximal nest curve $\C A$ into a pair of connected
components permuted by the complex conjugation, and thus, $\C
A\stmin\d M$ is connected, provided $\d M\neqqq\R A$, which is
the case for $d\ge 5$.
So, it is only left to check
that the group $\pi_1(\Cp2\stmin(\C A\cup M))$ is abelian.

There are several ways to check it.
For instance, one can refer to my old work \cite{F1} containing
computation of the homotopy type of $\Cp2\stmin(\C A\cup\Rp2)$
 and, in particular, of its fundamental group
(see also \S4 in \cite{FKV}).
This computation concerns a real curve
$\C A\subset\Cp2$ if it is an {\it $L$-curve}, i.e.,
$\C A$ can be obtained by a non-singular
perturbation from a curve $\C A_0=\C L_1\cup\dots\C L_d$
splitting into $d$ real lines, $\C L_i$, in a generic position.
The maximal nest curves, $\C A\subset\Cp2$, can be easily constructed
as $L$-curves, and the result of \cite{F1} gives a
presentation
$\pi=\pi_1(\Cp2\stmin(\C A\cup\Rp2))=\la a,b\,|\,a^{d}b^{d}=1\ra$,
where $a$, $b$ are represented by loops around
the two connected components of $\C A\stmin\R A$.
More specifically,
a basis point and these loops can be taken on the conic
$C=\{x^2+y^2+z^2=0\}\subset\Cp2$, which have the real point set empty.
The group $\pi_1(\Cp2\stmin(\C A\cup M))$ is obtained from $\pi$
by adding the relations corresponding to
puncturing the components $R_i$, $0\le i\le k$, $i\ne 1$,
of $\Rp2\stmin\R A$ (here $d=2k$ or $d=2k+1$).
Such a relation (as we puncture $R_i$) is $a^{d-i}b^i=b^{d-i}a^i=1$,
see \cite{F1}, or \S4 in \cite{FKV}.
A pair of the relations for $i=2$ and $i=3$ implies that $a=b$.

The arguments from \cite{F1} and \cite{FKV} 
 relevant
to the above calculation are briefly summarized in the Appendix.
\qed

\begin{rem}
It follows from the proof above that $\pi_1(\Cp2\stmin(\C A\cup M))$
is not abelian and $\C A\stmin\d M$ is not connected
for a maximal nest quartic, $\C A$.
\end{rem}

\section{The double surgery in the double covering}

\noindent{\bf Proof of Proposition \ref{1.4}.}
 The proof is based on the following two observations.
First, we notice that $Y(K,\phi)$ is obtained from $Y$ by a pair
of FS-surgeries along the tori parallel to $T$, then we notice
that such pair of surgeries is equivalent to a single FS-surgery
along $T$. The both observations are corollaries of Lemma 2.1 in
\cite{FS2}, so, I have to recall first the construction from
\cite{FS1}, \cite{FS2}.

An FS-surgery \cite{FS1} on a $4$-manifold $X$ along a torus
$T\subset X$, with the self-intersection
$T\circ T=0$, via a knot $K\subset S^3$ is defined as
a {\it fiber sum} $X\#_{T= S^1\times m_K} S^1\times M_K$,
that is an amalgamated connected sum of $X$ and $S^1\times M_K$
along the tori $T$ and $S^1\times m_K\subset  S^1\times M_K$.
Here $M_K$ is a $3$-manifold obtained by the $0$-surgery along $K$ in $S^3$,
and $m_K$ denotes a meridian of $K$ (which may be seen both in
$S^3$ and in $M_K$).
 Such a fiber sum operation can be viewed as a direct product
of $S^1$ and the corresponding $3$-dimensional operation,
which I call {\it $S^1$-fiber sum}.

More precisely, $S^1$-fiber sum
$X\#_{K=L}Y$ of oriented $3$-manifolds
$X$ and $Y$ along oriented framed knots $K\subset X$ and $L\subset Y$
 is the manifold obtained by gluing the complements
$\Cl(X\stmin N(K))$ and $\Cl(Y\stmin N(L))$ of tubular neighborhoods,
$N(K)$, $N(L)$, of $K$ and $L$ via a diffeomorphism
$f\:\d N(K)\to\d N(L)$ which identifies the longitudes
of $K$ with the longitudes of $L$ preserving their orientations,
and the meridians of $K$ with the meridians of $L$ reversing the
orientations.
 As it is shown in  Lemma 2.1 of \cite{FS2}, tying a knot $K$ in an arc
in $D^3$ can be interpreted as a fiber sum
$D^3\#_{m=m_K}M_K$, where $m$ is a meridian around this arc.
The meridians $m$ and $m_K$
are endowed here with the $0$-framings ($0$-framing of a meridian
makes sense as a meridian lies in a small $3$-disc).
 To understand this observation, it is useful to view
an $S^1$-fiber sum with $M_K$
as surgering a tubular neighborhood, $N(m)$, of $m$
and replacing it by the complement,
 $S^3\stmin N(K)$ of a tubular
neighborhood, $N(K)$, of $K$, so that the longitudes of $m$ are
glued to the meridians of $K$ and the meridians of $m$ to the
longitudes of $K$.
The framing of an arc in $D^3$ is preserved under such a fiber  sum,
so tying a knot in
the band $b\subset D^3$ is equivalent to taking
an $S^1$-fiber sum with $M_K$ along
a meridian $m_b$ around $b$.

The double covering over $D^3$ branched along $f$
is a solid torus, $N\cong S^1\times D^2$, and the pull back of $m_b$
splits into a pair of circles, $m_1,m_2\subset N$, parallel to
$m=S^1\times\{0\}$.
Therefore,
$Y(K,\phi)$ is obtained from $Y$ by performing FS-surgery twice,
along the tori
$$
T_i=S^1\times m_i\subset p^{-1}(U)\cong S^1\times N,\ \
i=1,2
$$
The following Lemma implies that this gives the same result as
a single FS-surgery along $T=p^{-1}(M)$
via the knot $K\#K$.
\qed

\begin{lem}\label{double-fiber-sum}
For any pair of knots, $K_1, K_2$, the manifold
$$M_{K_1}\#_{m_{K_1}=m_1}
N\#_{m_2=m_{K_2}}M_{K_2}$$
obtained by taking an $S^1$-fiber sum twice,
is diffeomorphic to
$N\#_{m=m_{K}}M_{K}$, for  $K=K_1\# K_2$,
via a diffeomorphism identical on $\d N$.
\end{lem}

\begin{proof}
A solid torus $N$ can be viewed as the complement $N=S^3-N'$ of an
open tubular neighborhood $N'$ of an unknot, so that $m,m_1,m_2$ represent
meridians of this unknot. Taking a fiber sum of $S^3$ with
 $M_{K_i}$ along $m_i=m_{K_i}$
is equivalent to knotting  $N'$ in $S^3$ via $K_i$.
So, performing
 $S^1$-fiber sum twice, along $m_1$ and $m_2$, we obtain the same result
as after taking fiber sum along $m$ once, via $K=K_1\# K_2$.
\end{proof}

\begin{rem}
The above additivity property can be equivalently stated as
$$
M_{K_1}\#_{m_{K_1}=m_{K_2}}M_{K_2}\cong M_{K_1\#K_2}
$$
\end{rem}

\noindent{\bf Proof of Proposition \ref{c-embedded}.}
 Lemma \ref{fg-cyclic} implies that, in the assumptions of
Proposition \ref{1.2},
$\pi_1(Y\stmin(F\cup T))$ is a cyclic group with
 a generator represented by a loop around $F$.
Thus, $\pi_1(Y\stmin T)=0$ and, by the Alexander duality,
$ H_3(Y,T)=H^1(Y\stmin T)=0$, which implies that
$[T]\in H_2(Y)$ has infinite order.

\begin{figure}[htb]
\includegraphics{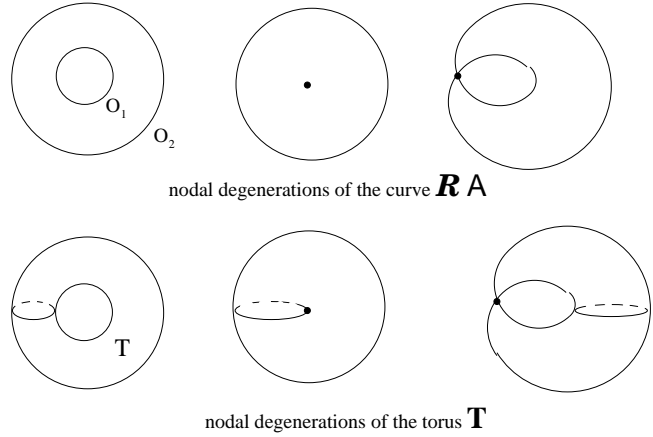}
\caption{Nodal degenerations of $\R A$ providing
$(-1)$-framed $D^2$-membranes on $T$
 $f_K$ into an unknotted arc}
\label{nodal}
\end{figure}

To check that $T$ is c-embedded it is enough to observe that
there exists a pair of vanishing cycles on $T$,
or more precisely, a pair of
$D^2$-membranes, $D_1,D_2\subset Y$, on $T$,
having $(-1)$-framing and intersecting at
a unique point $x\in T$, so that
$[\d D_1],[\d D_2]$ form a basis of $ H_1(T)$.
In the setting of Proposition \ref{1.3},
$Y\to\Cp2$ is a double covering branched along a maximal nest curve
$\C A$ and $T$ is a connected component of the real part of $Y$
(with respect to a certain real structure on $Y$ lifted from $\Cp2$).
Two nodal degenerations of $\C A$ shown on the top part of Figure \ref{nodal}
give nodal degenerations of the double covering $Y$.

In the first of the degenerations of $\C A$, a node appears as an
oval $O_1$ is collapsed into a point. In the second degeneration a
crossing-like node can be seen as the fusion point of the ovals
$O_1$ and $O_2$. Existence of such degenerations for our
explicitly constructed curve $\C A$ is known and trivial. Another
simple observation (which is obvious for quartics and thus follows for
any maximal nest curve of a higher degree)
is that our pair of nodal degenerations can be
united into one cuspidal degeneration. This means in particular
that the two vanishing cycles in $Y$
intersect transversally at a single point.

Furthermore, our complex vanishing cycles in $Y$
can be chosen $\conj$-invariant. Being a $(-2)$-sphere, each of
such complex cycles is divided by its real pair into a pair
of $(-1)$-discs. Choosing one disc from each pair, we obtain $D_1$
and $D_2$ that we need.

It is easy to view these $(-2)$-spheres and the
$(-1)$-disks explicitly.
First, note that $R_0$ is a $(-1)$-membrane on $\C A$ and
$p^{-1}(R_0)$ is the first of the $\conj$-symmetric vanishing cycles.
The $(-1)$-disk $D_1$ is any of its halves.
Furthermore, there is another $(-1)$-disk membrane, $Q$ on $\C A$
corresponding to the second nodal degeneration.
It can be chosen $\conj$-invariant and then
is split by $Q\cap\Rp2$ into semi-discs $Q=Q_1\cup Q_2$
 permuted by $\conj$.
$Q_i$ is bounded by the arcs $Q\cap\Rp2$ and $Q_i\cap\C A$.
The disk $D_2$ is any of the discs $p^{-1}(Q_i)$.
\qed

\section
{The case of $d$-fold branched covering}

Consider as before  a maximal nest curve, $\C A\subset\Cp2$,
of degree $d\ge2$, and $\C A_{K,\phi}$
obtained from $\C A$ via an annulus rim-surgery along $R_1$,
but now let us denote by $p\:Y\to\Cp2$ and $Y(K,\phi)\to\Cp2$
the {\it $d$-fold coverings}
branched along $\C A$ and $\C A_{K,\phi}$ respectively.
Consider a $d$-fold covering $N\to D^3$ branched along $f$.
The pull-back of $m_b$ consists of $d$ circles, $m_1,\dots,m_d$,
which are cyclically ordered.
Using a homeomorphism
$(D^3,f)\cong(D^2\times [0,1],\{z_1,z_2\}\times [0,1])$, where
$\{z_1,z_2\}\subset\Int(D^2)$, we present
$N$ as $F\times[0,1]$, where $F$
is a sphere with $d$ holes. The circles $m_i$ go around these holes.
An annulus rim-surgery in $\Cp2$ along $m_b\times S^1\subset D^3\times S^1$,
is covered by $d$ copies of FS-surgery along the tori
$T_i=m_i\times S^1\subset N\times S^1$.

The following observation implies that the Fintushel-Stern formula
for Seiberg-Witten invariants can be applied in this setting.

\begin{prop}\label{d-fold-tori}
Each of the tori $T_i$ is primitively c-embedded in the complement
of the others.
\end{prop}

\begin{proof}
A pair of $(-1)$-disc membranes, $D_1^i$, $D_2^i$,
on each of $T_i$ is constructed like
in the proof of Proposition \ref{c-embedded}.
Namely, $p^{-1}(R_0)$ consists of $d$ disks which yield the disks
$D_1^i$, that are
glued along $\{\text{pt}\}\times S^1\subset m_i\times S^1$.

Furthermore, $p^{-1}(Q_1)$ splits also into $d$ disks,
$Q_1^1,\dots,Q_1^d$. Let us choose their orientations induced from
a fixed orientation
of $Q_1$ and cyclically order in accord with the ordering of $T_i$, then
the unions $Q_1^i\cup(-Q_1^{i+1})$ provide the required discs $D_2^i$,
which are glued along $m_i\times \{\text{pt}\}$.
More precisely, $D_1^i$ are the parts of the components of
$p^{-1}(R_0)$ bounded by the intersections of the components
with the tori $T_i$, whereas $D_1^i$ are obtained from
 $Q_1^i\cup(-Q_1^{i+1})$ by a small shift making them membranes on $T_i$.
\end{proof}

Next, we observe that there exists only one linear dependence relation
between the classes $[T_i]\in H_2(Y)$.

\begin{prop}\label{lin-dependence}
The inclusion map $H_2(\bigcup_i T_i)\to H_2(Y)$ has kernel $\Z$
generated by the relation $\Sigma_{i=1}^d[T_i]=0$.
Here $T_i$ are oriented uniformly in accord with some fixed
orientation of $m_b\times S^1$.
\end{prop}

\begin{proof}
It is enough to show that $\pi_1(Y\stmin(N\times S^1))=0$, since
it implies that $H_3(Y,N\times S^1)\cong H^1(Y\stmin(N\times
S^1))=0$ and thus the inclusion map $H_2(N\times S^1)\to H_2(Y)$
is monomorphic. The first inclusion map in the composition
$H_2(\bigcup_i T_i)\to H_2(N\times S^1)\to H_2(Y)$ that we
analyze, is just $H_1(\d F)\otimes H_1(S^1)\to H_1(F)\otimes
H_1(S^1)$, and has kernel $H_2(F,\d F)\otimes H_1(S^1)\cong\Z$,
as stated in the Proposition.

Now note that $p^{-1}(R_1)$ is a deformational retract (spine)
of $N\times S^1$, so it is enough
to check the triviality of $\pi_1(Y\stmin(p^{-1}(R_1))$. This triviality
follows from that $\pi_1(\Cp2\stmin(\C A \cup R_1))$ is $\Z/d$,
 with a generator represented by a loop around $\C A$
(say, by the computation in \cite{F1} reproduced in
the Appendix),
and thus $\pi_1(Y\stmin p^{-1}(\C A \cup R_1))=0$.
\end{proof}

Proposition \ref{lin-dependence} together with the Fintushel-Stern formula
\cite{FS1}
guarantees that the Seiberg-Witten invariants of $Y(K,\phi)$ are distinct
for some sequence of knots $K$ with increasing degrees of $\Delta_K(t)$.

\noindent{\bf Proof of Theorem \ref{general}}
The case of a primitive class $[F]\in H_2(X)$ is considered in \cite{FS2}.
More precisely, the assumptions in Theorem 1.1
 in \cite{FS2} are satisfied
because our condition on the fundamental group yields that
$\pi_1(X\stmin F)$ is abelian and thus trivial,
existence of an irreducible
deformation of $F$ implies that $F\circ F\ge0$,
and $X_9$-degeneration guarantees that $F$ is not a rational curve.

If $[F]$ is divisible by $d\ge2$, then we consider a $d$-fold covering,
$p\:Y\to X$, branched along $F$ and perform
an annulus rim-surgery on $F$ along a membrane $M$ defined as follows.
Consider a local topological model of the singularity $X_9$,
 defined in $\C^2$ by the equation $(x^2+y^2)(x^2+2y^2)=0$, and a model
of its perturbation,
$(x^2+y^2-4\e)(x^2+2y^2-\e)=\delta$, where
$\e,\delta\in\R$, $0<\!\!<\delta<\!\!<\e<\!\!<1$.
The real locus of a perturbed singularity contains a pair of ovals
which bound together in $\R^2$ an annulus that we take as $M$.

The assumptions of Theorem \ref{general} imply those of Proposition
\ref{1.2}.
Namely, irreducibility of $F_0$ implies that $F\stmin {\d M}$ is connected
and commutativity of $\pi_1(X\stmin F_0)$ implies commutativity of
$\pi_1(X\stmin (F\cup M))$ via Van Kampen theorem.
Moreover, the singularity $X_9$ provides the topological picture
that was used in the above proof of Theorem \ref{thm-for-plane},
in the case of
$d$-fold covering. Namely, $X_9$ yields the both $(-1)$-disk membranes that
were used to show that the Fintushel and Stern formula can be applied to
$Y$.
\qed

\begin{rem}
Note that to apply the formula \cite{FS2} it is not required that
$b_2^+(Y)>1$. Nevertheless, it is so, because $b_2^+(Y)\ge d$, which
can be proved by observing
$d$ linearly independent pairwise orthogonal classes in $H_2(Y)$,
having non-negative squares. One of these classes is $[F]$, and the other
$(d-1)$ come from $p^{-1}(M)$, due to Proposition \ref{lin-dependence}
(each of these $(d-1)$ classes has self-intersection $0$).
\end{rem}

\section{Appendix: The topology of
$\Cp2\stmin(\Rp2\cup\C A)$ for $L$-curves $\C A$}
{\small

Let $\C A_0=\C L_1\cup\dots\cup\C L_d\subset\Cp2$ denote
the complex point set of a real curve of
degree $d$ splitting into $d$ lines, $\C L_i$.
Put $\widetilde V=C\cap\C A_0$, where $C$ is the conic from
the proof of Proposition \ref{1.3}.
 Our  first observation is that $C\stmin \widetilde V$ is
 a deformational retract of
 $\Cp2\stmin(\Rp2\cup\C A_0)$,
and moreover, the latter complement is homeomorphic to
$(C\stmin \widetilde V)\times \Int(D^2)$.
 To see it, it suffices to note that $\Cp2\stmin\Rp2$ is fibered
 over $C$ with a $2$-disc fiber, each fiber being
a real semi-line, that is a connected component of $\C L\stmin\R
L$ for some real line $\C L\subset\Cp2$, where $\R L=\C
L\cap\Rp2$. This fibering maps a semi-line into its intersection
point with $C$.

It is convenient to view
the quotient $C/\conj$ of the conic $C$ by the complex conjugation
as the projective plane, $\widehat{\Rp{}}^2$,
dual to $\Rp2\subset\Cp2$, since each real line, $\C L$,
intersects $C$ in a pair of conjugated points.
If we let
 $V=\{l_1,\dots,l_d\}\subset\widehat{\Rp{}}^2 $ denote the set of points
$l_i$ dual to the lines $\R L_i\subset\Rp2$, then
$\widetilde V=q^{-1}(V)$, where
$q\: C\to C/\conj$ is the quotient map.

The information about a perturbation of $\C A_0$ is encoded in a
{\em genetic graph of a perturbation},
$\G\subset\widehat{\Rp{}}^2$. The graph $\G$ is a complete graph
with the vertex set $V$, whose edges are  line segments. Note that
there exist two topologically distinct perturbations of a real
node of $\R A_0$ at $p_{ij}=\R L_i\cap\R L_j$, as well as there
exist two line segments in $\widehat{\Rp{}}^2$ connecting the vertices
$l_i, l_j\in V$. Let $\R A$ denotes a real curve obtained from $\R
A_0$ by a sufficiently small perturbation. Then the edge of $\G$
connecting $l_i$ and $l_j$ contains the points dual to those lines
passing through $p_{i,j}$ which do not intersect $\R A$ locally,
in a small neighborhood of $p_{i,j}$.

\begin{figure}[htb]
\includegraphics{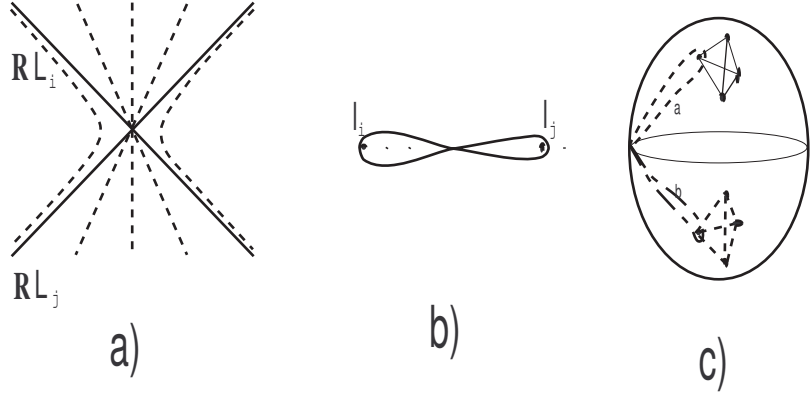}
\caption{a) A perturbation of a real node; the dashed lines are dual
to the points of an edge of $\G$;
b) A figure-eight loop along an edge of $\G$;
c) The loops in $C\stmin \widetilde V$ representing
generators ``a'' and ``b''}
\label{l-curves}
\end{figure}

The complement $\Cp2\stmin(\C A\cup\Rp2)$ turns out to be
homotopy equivalent to a $2$-complex obtained from
$C\stmin \widetilde V$ by adding $2$-cells glued along a figure-eight
shaped loops along the edges of $\widetilde{\G}=q^{-1}(\G)\subset C$.
Such $2$-cells identify pairwise certain
generators of $\pi_1(C\stmin\widetilde V)$
``along the edges'' of $\widetilde{\G}$ (cf. \cite{FKV} for details).
This easily implies that the group
$\pi_1(\Cp2\stmin(\C A\cup\Rp2))$ is generated by a pair of elements,
$a$ and $b$, represented by a pair of loops in $C\stmin \widetilde V$ around
a pair of conjugated vertices of $\widetilde V$.

For example, for a maximal nest curve, the graph $\G$ is contained
in an affine part of $\widehat{\Rp{}}^2$, i.e., has no common points
with some line in $\widehat{\Rp{}}^2$, namely, with a line
dual to a point inside the inner oval of the nest.
Therefore, the graph $\widetilde{\G}$ splits into two connected components
separated by a big circle in $C$.
 A loop around any vertex of $\widetilde V$ from one of these components
 represents $a$, and a loop around a vertex from the other
 component represents $b$.
It is trivial to observe also the relation $a^db^d=1$
(which is indeed a unique relation in the case of maximal nest curves).

As we puncture $\Rp2$ at a point $x\in \Rp2\stmin\R A_0$,
we attach a $2$-cell to $C\stmin \widetilde V$ along the big circle
$S_x\subset C$ dual to $x$.
If $x$ moves across a line $\R L_i$, then $S_x$ moves across the pair
of points $q^{-1}(l_i)$.
 Since a small perturbation and puncturing are located at
distinct points of $\Cp2$ and can be done independently,
it is not difficult to see that if we choose $x\in R_i$
(in the case of a maximal nest curve $\C A$), then the big circle $S_x$
cuts $C$ into the hemispheres, one of which
contains $i$ vertices from one component of $\widetilde{\G}$ and
$d-i$ vertices from the other component.
This gives relations $a^ib^{d-i}=a^{d-i}b^i=1$.
}

\end{document}